\documentclass[12pt]{amsart}
\usepackage{verbatim}
\usepackage{eucal}
\usepackage{amssymb}
\usepackage{mathrsfs}
\usepackage{graphicx}
\usepackage{psfrag}
\newif \ifwide
\newif \ifavnermargin
\def \makemargins{
\ifwide
        \oddsidemargin .25in
        \evensidemargin .25in
        \textwidth 6.00in
\else
\fi
\ifavnermargin
        \headheight=7pt
        \textheight=574pt
        \textwidth=432pt
        \topmargin=14pt
        \oddsidemargin=18pt
        \evensidemargin=18pt
\else
\fi
}

\widetrue
\makemargins

\theoremstyle{plain}
\newtheorem{theorem}[subsection]{Theorem}
\newtheorem{proposition}[subsection]{Proposition}
\newtheorem{lemma}[subsection]{Lemma}
\newtheorem{corollary}[subsection]{Corollary}

\theoremstyle{definition}
\newtheorem{definition}[subsection]{Definition}
\newtheorem{example}[subsection]{Example}

\theoremstyle{remark}
\newtheorem{remark}[subsection]{Remark}

\newcount\timehh\newcount\timemm
\timehh=\time
\divide\timehh by 60 \timemm=\time
\newcommand{\draftauthor}[1]{\author{#1
    {
      --- \protect \protect\sc\today\ ---
      \ifnum\timehh<10 0\fi\number\timehh\,:\,\ifnum\timemm<10 0\fi\number\timemm
      \protect \, \, \protect \bf DRAFT
    }
  }
}
\newcommand{\RR}{{\mathbb R}}
\newcommand{\CC}{{\mathbb C}}
\newcommand{\ZZ}{{\mathbb Z}}

\newcommand{\PP}{{\mathbb P}}

\newcommand{\ee}{{\mathrm{e}}}
\newcommand{\ii}{{\mathrm{i}}}

\newcommand{\HHH}{{\mathfrak{H}}}
\newcommand{\TTT}{{\mathscr{T}}}

\newcommand{\EEE}{{\mathscr{E}}}
\newcommand{\MMM}{{\mathscr{M}}}
\newcommand{\SSS}{{\mathscr{S}}}
\newcommand{\QQQ}{{\mathscr{Q}}}

\newcommand{\stroke}[2]{{#1\!\bigm|\!{#2}}}
\renewcommand{\theta}{\vartheta}
\renewcommand{\tilde}{\widetilde}

\renewcommand{\mod}{\bmod\,}

\DeclareMathOperator{\Hom}{Hom}

\DeclareMathOperator{\Sum}{Sum}
\DeclareMathOperator{\cone}{cone}
\begin{document}

\title{Toric modular forms of higher weight}

\newif \ifdraft
\def \makeauthor{
\ifdraft
        \draftauthor{Lev A. Borisov and Paul E. Gunnells}
\else
\author{Lev A. Borisov}
\address{Department of Mathematics\\
Columbia University\\
New York, NY  10027}
\email{lborisov@math.columbia.edu}

\author{Paul E. Gunnells}
\address{Department of Mathematics and Statistics\\
University of Massachusetts\\
Amherst, MA  01003}
\email{gunnells@math.umass.edu}

\fi
}

\draftfalse
\makeauthor

\ifdraft
        \date{\today}
\else
        \date{March 16, 2002}
\fi


%
%

\begin{abstract}
In the papers \cite{vanish, toric} we used the geometry of
complete polyhedral fans to construct a subring $\TTT (l)$ of the modular
forms on $\Gamma _{1} (l)$, and showed that for weight two the
cuspidal part of $\TTT (l)$ coincides with the space of cusp forms of
analytic rank zero.  In this paper we show that in weights greater
than two, the cuspidal part of $\TTT (l)$ coincides with the space of
all cusp forms.
\end{abstract}
\maketitle

%
%
\section{Introduction}\label{introduction}
\subsection{}
In \cite{vanish, toric} we used the geometry of complete polyhedral
fans to construct a subring $\TTT_* (l)$ of the modular forms on
$\Gamma _{1} (l)$.  If $l\geq 5$, we showed that $\TTT_* (l)$ is
generated in weight one by certain Eisenstein series, and in
\cite[Theorem 4.11]{vanish} we showed that for weight two the cuspidal
part of $\TTT_* (l)$ coincides with the space of cusp forms of
analytic rank zero.  The main result of this paper, Theorem
\ref{main}, is that in weights greater than two, the cuspidal part of
$\TTT_* (l)$ coincides with the space of all cusp forms.  In fact, we
prove a stronger statement: we define certain weight $k$ toric modular
forms $\tilde s^{(k)}_{a/l}$, and show that any cusp form can be
written as a $\CC$-linear combination of the forms $\tilde s^{(k)}_{a/l}$ and
pairwise products of the form $\tilde s^{(m)}_{a/l} \tilde{s}_{b/l}^{(n)}$, 
where $m+n=k$ and $m,n>0$.

The proof of Theorem \ref{main} is formally very similar to the proof
of \cite[Theorem 4.11]{vanish}.  Let $\SSS (l)$ be the space of weight
$k$ holomorphic cusp forms on $\Gamma_{1} (l)$.  We define a map $\rho
\colon \SSS (l)\rightarrow \SSS (l)$, and show that its image contains
all newforms for $k\geq 3$. We describe the map $\rho$ in terms of
Manin symbols, which allows us to write $\rho(f)$ in terms of products
of certain explicit toric Eisenstein series.  A key role is played by
certain weight $k$ Manin symbols $\{R_{(m,n)} | m,n\in \ZZ\}$ that
satisfy relations similar to weight two Manin symbols.

\subsection{}
Here is an outline of the paper.  In Section \ref{s1} we review
results about toric modular forms, and in Section \ref{s2} review results
about Manin symbols and introduce the symbols $R_{(m,n)}$.  In Section
\ref{s3} we describe \emph{$(\mod l)$-polynomials}, a technical tool
we use later to manipulate $q$-expansions.  We prove the main result
along with some corollaries in Section \ref{s4}.

The remaining sections contain complements to the main result and
results proved in \cite{vanish}.  In Section \ref{s.toricmap} we use
products of Eisenstein series of higher weight to define a map $\mu$
from weight $k$ Manin symbols to a certain quotient of the space of
weight $k$ modular forms.  This map is analogous to the map $\mu$ in
\cite[Definition 3.11]{vanish}, but some complications do occur in the
higher case.  Finally, in Section \ref{s.Hecke} we show that the map
from symbols to forms is compatible with the action of the Hecke
operators.

Throughout the paper we keep our arguments as elementary as
possible. In particular, we avoid using results of \cite{toric} that
are based on the Hirzebruch-Riemann-Roch theorem for toric
varieties. While this complicates the proofs a bit, it makes the paper
accessible to readers with no knowledge of toric varieties. We outline
an alternative approach to the results using toric geometry in Remarks
\ref{toric1} and \ref{toric2}.

\section{Toric modular forms}\label{s1}
\subsection{}
We briefly review the definition of toric forms; more details can be
found in \cite{vanish, toric}.  Let $k$ and $l$ be positive integers,
and suppose $l\geq 5$.  As usual let $q=\ee^{2\pi \ii \tau}$, where
$\tau$ is in the upper halfplane $\HHH$.  A holomorphic modular form
of weight $k$ on the group $\Gamma _{1} (l)$ is called \emph{toric} if
it can be expressed as a homogeneous polynomial of degree $k$ in the
functions $\tilde s_{a/l}(q)$ given by
$$
\tilde s_{a/l}(q):=(\frac12-\frac al)+\sum_{n>0}q^n\sum_{d|n}
(\delta_d^{a\mod l}-\delta_d^{-a\mod l}),\quad a=1,\dots ,l-1.
$$
Here $\delta_d^{a\mod l}$ is $1$ if $a=d\mod l$ and is $0$ otherwise.
The space of toric forms $\TTT_{*}(l)$ of all weights is thus
generated as a graded ring by certain weight one Eisenstein series.
By results of \cite{vanish, toric}, $\TTT_{*}(l)$ is known to be stable
under the action of the Hecke operators, and under Atkin-Lehner
lifting.

\begin{proposition}\label{unnamed}
The ring $\TTT_{*} (l)$ contains the modular forms
\begin{equation}\label{stilde}
\tilde s_{a/l}^{(k)} := C + \sum_{n>0}q^n\sum_{d|n}d^{k-1}
(\delta_d^{a\mod l}+(-1)^k\delta_d^{-a\mod l}),
\end{equation}
where
$k\geq 2$ and $a=0,\dots ,l$,
except for $(a,k) = (0,2)$. Here $C$ is a constant determined uniquely
by the modularity of $\tilde s_{a/l}^{(k)}$.
\end{proposition}

\begin{proof}
Let $\theta (z,\tau)$ be Jacobi's theta function \cite{Chandra}.  Then
it is easy to construct a given $\tilde s_{a/l}^{(k)}$
as a linear combination of the modular forms
\begin{align*}
s_{a/l}^{(k)}(q)&:= (2\pi\ii)^{-k}
\left(\frac{\partial^k}{\partial z^k}\right)_{z=0}{\log}
\left(
{
\frac{z\theta(z+a/l,\tau)\theta'(0,\tau)}
{\theta(z,\tau)\theta(a/l,\tau)}
}
\right)\\&=
C+\sum_{n>0}q^n\sum_{d|n}d^{k-1}(\ee^{2\pi\ii da/l} + (-1)^{k}\ee^{-2\pi\ii
da/l} - 2\delta_k^{0\mod 2})
\end{align*}
of \cite[Section 4.4]{toric}, and the standard level one Eisenstein series
\begin{equation}\label{eis.series}
E_{k} := C + \sum _{n>0}q^n\sum_{d|n}d^{k-1}
\end{equation}
for even $k\geq 4$.  The forms $s_{a/l}^{(k)}$ and $E_{k}$ for $k>2$
are toric, see \cite[Theorem 4.11 and Remark 4.13]{toric}.  We use
here and in what follows the convention denoting constant terms whose
exact value is irrelevant by $C$.
\end{proof}

\subsection{}
We must exclude $(a,k) = (0,2)$ from the statement since $E_2$ is not
modular.  However, it will be convenient to allow $\tilde
s_{0/l}^{(2)}$ in later arguments, which merely amounts to working in
the larger ring $\TTT_*(l)[E_2]$.  In fact, since we will never multiply
more than two of the $\tilde s$'s together, we will be working in
$\TTT_*(l)+E_2\TTT_*(l)$. We call elements of this ring \emph{toric
quasimodular forms} (cf. \cite{Zagier.elliptic}).

\begin{remark}\label{lessthanfive}
The statement of Proposition \ref{unnamed} is true for all $l=2,3,4$, but the
definition of toric modular forms there is a bit more complicated. However,
it turns out that for these levels all modular forms are toric, as defined
in \cite{toric}. From now on we will call any polynomial in 
$\tilde s^{(k)}_{a/l}$ a toric quasimodular form of level $l$. 
Then the rest of the paper works for an arbitrary level $l\geq 1$.
\end{remark}

\begin{definition}\label{pairs}
By a slight abuse of notations 
we say that a weight $k$ quasimodular form $f$ can be written as a
\emph{linear combination of pairs} if $f$ can be written as a
$\CC$-linear combination of the forms $\tilde s_{a/l}^{(k)}$ and 
$\tilde s_{a/l}^{(m)} \tilde s_{b/l}^{(n)}$ 
where $m,n>0$, $m+n=k$, and $a,b=0,\dotsc ,l-1$.
\end{definition}

\begin{proposition}\label{derivs.are.there}
The space of toric quasimodular forms contains the derivatives
$\partial_{\tau} \tilde s^{(k)}_{a/l}$.  Moreover, each $\partial_{\tau}
\tilde s^{(k)}_{a/l}$ can be written as a linear combination of pairs.
\end{proposition}

\begin{proof}
The span of $\tilde s^{(k)}$ is the same as the span of $s^{(k)}$ and $E_{k}$,
so we will instead consider their derivatives.
The $q$-expansion of $\partial_{\tau} s^{(k)}_{a/l}$ is
\begin{equation}\label{deriv.q.exp}
2\pi \ii \sum_{n>0}q^n\sum_{d|n}nd^{k-1}(\ee^{2\pi\ii da/l} + (-1)^{k}\ee^{-2\pi\ii
da/l} - 2\delta_k^{0\mod 2}).
\end{equation}
Now take (cf. \cite[proof of Prop. 3.8]{vanish})
\begin{equation}\label{identity}
s_{\alpha}^{(2)}+s_{\alpha}^2 =
\frac 16 - 2\sum_{n>0}q^n\sum_{d|n}\frac nd(\ee^{2\pi\ii \alpha d}+
\ee^{-2\pi\ii \alpha d}),
\end{equation}
and differentiate it $k$ times with respect to $\alpha$.  Let
$F_{\alpha} (q)$ be the resulting right hand side.  It is easy to
express \eqref{deriv.q.exp} as a linear combination of $\{F_{a/l} (q)
\}$ and the derivatives $\partial_{\tau} E_{k}$, so it suffices to
show that these forms can be expressed as a linear combination of pairs.  

We consider first the derivatives
of $s_{\alpha}^{(r)}$ with respect to
$\alpha$.  Putting $E_{\rm odd} = 0$, we have
$$
(2\pi\ii)^{-1}\frac{\partial}{\partial \alpha} s_{\alpha}^{(r)}=
s_{\alpha}^{(r+1)} -
(2\pi\ii)^{-r-1}
\left(\frac{\partial^{r+1}}{\partial z^{r+1}}\right)_{z=0}\log
\left(
{
\frac{z\theta'(0,\tau)}{\theta(z,\tau)}
}
\right)
=s_{\alpha}^{(r+1)} - 2E_{r+1}, 
$$
and the statement follows from the fact that $E_r$ and $s^{(r)}$
can be written as linear combinations of $\tilde s^{(r)}$.  
For $\partial_{\tau} E_{k}$ we argue as follows. Expand both
sides of the equation \eqref{identity} in a Laurent series in 
$\alpha$ around $\alpha=0$. 
The coefficient at $\alpha^k$ on the right hand side
of \eqref{identity} is equal to $\partial_{\tau}E_k$, up to a 
multiplicative and an additive constant. To expand the left hand side
notice that, up to the terms constant in $q$,
the Laurent coefficient of $s_{\alpha}$ at $\alpha^k$ is 
a multiple of $E_{k+1}$, which follows from expanding
$\ee^{2\pi\ii d\alpha}$ in the definition of $s_{\alpha}$.
It is easy to see that $s_{\alpha}$ has a simple pole at $\alpha=0$
with a constant residue, so the coefficient of the Laurent
expansion of $s_{\alpha}^2$ at $\alpha^k$ is a linear combination
of $E_{r}E_{k+2-r}$, $E_{k+2}$ and some $E_r$ for $r<k+2$.
Then the modular transformation properties of $\partial_\tau E_k$
finish the argument.
\end{proof}

Next we describe the action of $\Gamma_0(l)/\Gamma_1(l)$ on $\tilde s$.
\begin{proposition}\label{gamma0action}
Let $\gamma\in \Gamma_0(l)$ have diagonal entries
$p^{-1}$ and $p$ $\mod l$ respectively. Then
$$\gamma \tilde s_{a/l}^{(k)}= \tilde s_{p^{-1}a/l}^{(k)}.$$
\end{proposition}

\begin{proof}
The transformation properties of $\theta$ (cf. \cite[Prop. 4.3]{toric}) imply 
$$\gamma s_{a/l}^{(k)}= s_{pa/l}^{(k)}.$$
One can then use linear combinations of the forms in the proof of
Proposition \ref{unnamed} to determine the action of $\gamma$ on $\tilde s$.
We leave the details to the reader.
\end{proof}

\section{Manin symbols}\label{s2}

\subsection{}\label{relationsection}
This section closely follows \cite{Merel}, to which the reader is
referred for more details.  Let $l>1$ be an integer, and let
$E_{l}\subset (\ZZ /l\ZZ )^{2}$ be the subset of pairs $(u,v)$ such
that $\ZZ u + \ZZ v = \ZZ /l\ZZ $.  The space of \emph{Manin symbols}
of weight $k$ and level $l$ is the $\CC$-vector space generated by the
symbols $x^ry^s(u,v)$, where $r$ and $s$ are nonnegative integers
summing to $k-2$ and $(u,v)\in E_l$, modulo the following relations:
\begin{enumerate}
\item $x^ry^s(u,v) + (-1)^rx^sy^r(v,-u) = 0$.
\item $x^ry^s(u,v) + (-1)^ry^r(x-y)^s(v,-u-v) +
         (-1)^s(y-x)^rx^s(-u-v,u) = 0$.
\end{enumerate}

We denote the space of Manin symbols by $M$ (we omit the level and
weight from the notation since it will be clear from the context).
Two subspaces of $M$ will play an important role in what follows.  Let
$\iota \colon M \rightarrow M$ be the involution 
$x^ry^s(u,v)\mapsto (-1)^rx^ry^s(-u,v)$.
\begin{definition}\label{plusnminus}
The space of \emph{plus symbols} $M_{+}\subset M$ is the subspace
consisting of symbols $w$ satisfying $\iota (w) = w$.  Similarly,
the space of \emph{minus symbols} $M_{-}\subset M$ is the subspace
consisting of symbols $w$ satisfying $\iota (w) = -w$.
\end{definition}

We have \emph{symmetrization maps} $(\phantom{a},\phantom{a})_{\pm }\colon M
\rightarrow M_{\pm }$ given by $x^ry^s(u,v)_{\pm } :=
(x^ry^s(u,v)\pm (-1)^rx^ry^s(-u,v))/2$.

\subsection{}\label{degenerate}
Let $M^{*} = \Hom_{\CC} (M,\CC)$ be the dual of the space of Manin
symbols. For any $\varphi \in M^{*}$, we define $\varphi $ on ``degenerate''
symbols $x^ry^s(u,v)$ with
$\ZZ u+\ZZ v \not = \ZZ /l\ZZ $ by setting $\varphi(x^ry^s(u,v)) = 0$.
This convention is somewhat artificial but turns out to be quite
useful.

\subsection{}\label{duality}
%
There exists a natural pairing between the spaces of Manin symbols
and the spaces of cusp forms, see \cite{Merel}.
Let $\MMM(l)$ be the $\CC$-vector space of weight $k$ holomorphic
modular forms on $\Gamma_{1} (l)$, and let $\SSS (l)\subset \MMM (l)$
be the subspace of cusp forms.  For $x^ry^s(u,v)\in M$ and $f \in \SSS (l)$ let
$g=\left(
\begin{smallmatrix}
a&b\\c&d
\end{smallmatrix}
\right)$ be an element of $\Gamma(1)$ with $(c,d)=(u,v)\mod l$. Then
the integral
\[
\int_{0}^{\ii\infty} (c\tau+d)^{-k}f(\frac{a\tau+b}{c\tau+d})\tau^r \,d\tau
\]
does not depend on the choice of $g$. Moreover, it is compatible
with the relations on modular symbols, and we obtain a pairing
\begin{align*}
M\times \SSS (l)&\longrightarrow \CC,\\
(x^ry^s(u,v), f)&\longmapsto \langle f, x^ry^s(u,v)\rangle.
\end{align*}
In general this pairing is degenerate, but one can identify a subspace
of cuspidal Manin symbols $S$ such that the pairing is non-degenerate
on $S_{\pm}\times \SSS(l)$.  We will not use this fact, but details
can be found in \cite{Merel}.

\subsection{}\label{hecke.action}
Next we present Merel's description of the Hecke action on the Manin
symbols.  Let $n\geq 1$ be an integer, and let $T_{n}$ be the
associated Hecke operator.  We denote the action of $T_{n}$ on a
modular form $f$ by $\stroke{f}{T_{n}}$.  For any positive integer
$n$, we define a set $H (n)\subset \ZZ^{4}$ by
\begin{equation}\label{merel.set.def}
H (n) = \{(a,b,c,d) \mid  a>b\geq 0, \quad 
d>c\geq 0,\quad 
ad-bc = n\}.
\end{equation}
\begin{theorem}\label{HeckeisHecke}
\cite[Theorem 2 and Proposition 10]{Merel} For any positive integer
$n$ coprime to $l$, define an operator $T'_{n} \colon M\rightarrow M$
by
\begin{equation}\label{heckeact}
T'_{n} \,x^ry^s(u,v) = \sum _{H (n)}(ax+by)^r(cx+dy)^s(au+cv, bu+dv).
\end{equation}
If $n$ is not coprime to $l$, then we define $T'_{n}$ by
\eqref{heckeact} but omit terms with  
$g.c.d.(l,au+cv,bu+dv ) > 1$.
Then $T'_{n}$ is the adjoint of $T_{n}$ with respect to the pairing
$\langle \phantom{a}, \phantom{a}\rangle$, that is
\[
\langle \stroke{f}{T_{n}},x^ry^s (u,v) \rangle= \langle f, T'_{n}\,
x^ry^s(u,v)\rangle.
\]
\end{theorem}
We will abuse notation in what follows and write $T_{n}$ for $T'_{n}$.
It is also proved in \cite{Merel} that this Hecke action is
compatible with the symmetrization maps:
\begin{proposition}\label{HeckeonSymbols}
We have
\[
T_{n} (x^ry^s(u,v)_{\pm }) = (T_{n}\,x^ry^s (u,v))_{\pm }.
\]
\end{proposition}

\subsection{}\label{rkls}
To conclude this section 
we associate to every pair of integers
$(m,n)$ a certain Manin symbol $R_{(m,n)}$. These symbols satisfy
relations analogous to those satisfied by the weight two symbols.

\begin{definition}\label{Rmn}
Let $m, n\in\ZZ $. If $g.c.d.(m,n,l)=1$ then we
let  $R_{(m,n)}=(mx+ny)^{k-2}(m,n)$.  If $g.c.d.(m,n,l)>1$ then we put $R_{(m,n)}=0$.
\end{definition}

We remark that even though $R_{(m,n)}$ is built out of the Manin
symbol $(m,n)$, its value depends on more than just the residues of $m$
and $n$ modulo $l$.  It is straightforward to see that the symbols
$R_{(m,n)}$ obey the following relations:
\begin{enumerate}
\item $R_{(m,n)}+R_{(-n,m)}=0$.
\item $R_{(m,n)}+R_{(-m-n,m)}+R_{(n,-m-n)} = 0$.
\item $\bigl(R_{(m,n)}\bigr)_{\pm}=(R_{(m,n)}\pm R_{(-m,n)})/2$.
\end{enumerate}
We denote the images of $R_{(m,n)}$ under the symmetrization maps by
$R_{(m,n)}^\pm$.

\section{(Mod $l$)-polynomials}\label{s3}

\subsection{} To simplify
later manipulations with $q$-expansions, we now introduce certain
functions.  Fix a positive integer $l$.

\begin{definition}
A function $h\colon \ZZ\to \CC$ is called a \emph{$(\mod l)$-polynomial}
if its restriction to each coset $l\ZZ + k$ is a polynomial.
\end{definition}

One can think of a $(\mod l)$-polynomial as a set of $l$
ordinary polynomials, one for each residue modulo $l$. For example, the
function that equals $m^2+m$ when $m$ is even and $m^3$ when $m$ is
odd is a $(\mod 2)$-polynomial.

The set of all $(\mod l)$-polynomials forms a ring.  One can
analogously define $(\mod l)$-polynomials $h(m,n)$ of two variables by
requiring polynomiality on each pair of cosets $(l\ZZ+k_1,l\ZZ+k_2)$.

\subsection{}
We say that a $(\mod l)$-polynomial $h$ is \emph{odd} if $h (-m) = -h
(m)$.  Note that the individual polynomials constituting an odd $(\mod
l)$-polynomial aren't independent, since the polynomials sitting over
the residues $a \mod l$ and $-a \mod l$ are related.  The space of odd $(\mod
l)$-polynomials will be of particular importance
to us, due to the following proposition.
\begin{proposition}\label{easyremark}
Let $h$ be an odd $(\mod l)$-polynomial. Then up to a constant, the function
$$
f(q) = \sum_{D>0}q^D\sum_{d|D}h(d)
$$
is a linear combination of $\{\tilde s_{a/l}^{(k)} \mid k\geq 1,\, a=1,\dotsc
l-1\}$.
Conversely, every linear combination of $\tilde s_{a/l}^{(k)}$ has the above
form, up to an additive constant.
\end{proposition}

\begin{proof}
Let $r_{a,k} (m)$ be the $(\mod l)$-polynomial given by
$$
r_{a,k}(m)=m^k\delta_m^{a\mod l}-(-1)^{k}m^k\delta_{m}^{-a\mod l}.
$$
Then any odd $(\mod l)$-polynomial is a linear combination of the
$r_{a,k}$'s.  The result follows easily
from the definition of $\tilde s_{a/l}^{(k)}$
in \eqref{stilde}.
\end{proof}

The following result allows one to construct odd one-variable $(\mod l)$-polynomials
from even two-variable $(\mod l)$-polynomials.
\begin{proposition}\label{conesum}
Let $G\colon\ZZ^2\to\CC$ be a two-variable $(\mod l)$-polynomial such that
$G(-n_1,-n_2)=G(n_1,n_2)$, and let $N$ be a positive integer. Then
$$
f(d):=\sum_{0<n<Nd}G(n,d)+\frac12 G(0,d)+\frac12 G(Nd,d)
$$
is an odd $(\mod l)$-polynomial.
\end{proposition}

\begin{proof}
First we note that the space of all even two-variable $(\mod
l)$-polynomials is spanned by the family of functions
$$
G(n_1,n_2) = n_1^rn_2^s \ee^{{2\pi\ii}(k_1n_1 + k_2 n_2)/l}
+(-n_1)^r(-n_2)^s\ee^{-{2\pi\ii}(k_1n_1 + k_2 n_2)/l}
$$
for nonnegative integers $r$ and $s$ and integers $k_1$ and $k_2$.
We use $\alpha_{i} = 2\pi \ii k_{i}/l$ to write such a $G$ as
$$
G(n_1,n_2)=n_1^rn_2^s(\ee^{\alpha_1n_1+\alpha_2 n_2}
+(-1)^{r+s}\ee^{-\alpha_1n_1-\alpha_2 n_2}).
$$
Now it suffices to treat the case $r=s=0$, since all others can then be
handled by partial differentiation with respect to $\alpha_1,
\alpha_2$.  For $r=s=0$ and $\ee^{\alpha_1}\neq 1$, an explicit
calculation gives
$$
f(d) = \sum_{0<n<Nd}(\ee^{\alpha_1n+\alpha_2d}+\ee^{-\alpha_1n-\alpha_2d})
+\frac 12 (\ee^{\alpha_2 d}+\ee^{-\alpha_2d})
+\frac 12 (\ee^{(\alpha_1N+\alpha_2) d}+\ee^{-(\alpha_1N+\alpha_2) d})
$$
$$
=\frac{(1+\ee^{\alpha_1})}{2(1-\ee^{\alpha_1})}
\Big(
\ee^{\alpha_2 d}-\ee^{(\alpha_1 N+\alpha_2)d}
-\ee^{-\alpha_2 d}+\ee^{-(\alpha_1 N+\alpha_2)d}
\Big).
$$
This is clearly an odd function in $d$, and after letting $\alpha_{i}
= 2\pi \ii k_{i}/l$ is obviously $(\mod l)$-polynomial in $d$.  The
case $\ee^{\alpha_1}=1$ follows by analytic continuation.
\end{proof}

\subsection{}
The following technical statement will be needed for the proof
of Lemma \ref{f2lemma}.

\begin{proposition}\label{oddprop}
Fix a weight $k$ cusp form $f$ on $\Gamma_1(l)$,
and define a function $h\colon \ZZ_{>0}\to \CC$ by 
$$h(m):=\langle f, R_{(m,0)}^+\rangle
+2\langle f, \sum_{m>i>0} R_{(m,m-i)}^+\rangle.$$
Then $h$ extends to an odd $(\mod l)$-polynomial. 
\end{proposition}

\begin{proof}
We use the symmetries of $R^+$ to rewrite $h(m)$ as
$$h (m)=\langle f, \sum_{-m<i<m} R_{(m,i)}^+\rangle=
\sum_{0<i<2m}\langle f, R_{(m,i-m)}^+\rangle
+\frac 12\langle f, R_{(m,-m)}^+\rangle
+\frac 12 \langle f, R_{(m,m)}^+\rangle.$$
Then Proposition \ref{conesum} finishes the proof.
\end{proof}

\section{Main theorem}\label{s4}

\subsection{}\label{rho}
Fix a weight $k\geq 3$ and a level $l$.  In this section we define
an endomorphism of the space $\SSS(l)$ of cusp forms of weight $k$
with respect to $\Gamma_1(l)$, and prove that its image contains all
newforms.  This definition is a generalization of
\cite[Definition 4.2]{vanish} to $k>2$.

\begin{definition}
Let
$\rho\colon \SSS(l)\to \SSS(l)$ be the linear map
$$\rho(f)=\sum_{n=1}^{\infty}\left(\int_{0}^{\ii\infty}(\stroke{f}{
T_n})(s)ds\right) q^n.$$
\end{definition}

\begin{proposition}\label{cuspform}
The form $\rho(f)$ is a cusp form with nebentypus equal to that of $f$.
\end{proposition}

\begin{proof}
The statement follows from \cite[Theorem 6]{Merel}; see also
\cite[Proposition 4.3]{vanish}.
\end{proof}

The map $\rho$ was used in \cite{vanish} because its image contains
all newforms of weight two whose $L$-functions don't vanish at the
center of the critical strip.  The analogous
statement for higher weights is the following:

\begin{proposition}\label{nothingvanishes}
The image of $\rho$ contains all newforms.
\end{proposition}

\begin{proof}
One needs to show that for any newform $f$
$$
\int_{0}^{\ii\infty}f(\tau)\,d\tau\neq 0,
$$
which is equivalent to $L(f,1)\neq 0$.  Without loss of generality we
may assume that $f$ is a Hecke eigenform.  If $k>3$ then $L (f,1)$ is a special
value outside the critical strip, and so cannot vanish by absolute
convergence of the Euler product.  If $k=3$ then $L (f,1)$ is a
special value on the boundary of the critical strip.  By \cite[Theorem
1.3]{jacquet} this special value cannot vanish.
\end{proof}

\subsection{}\label{keylemma}
Fix a cusp form $f$.  By Theorem \ref{HeckeisHecke} we can express $\rho(f)$ in terms of modular
symbols as
$$\rho(f)=\sum_{n=1}^{\infty}q^n\langle f, T_n y^{k-2}(0,1)\rangle
=\sum_{n=1}^{+\infty}q^n\langle f, T_n R_{(0,1)}^+\rangle
=\sum_{n=1}^{\infty}q^n\langle f, \sum _{H (n)}R_{(c,d)}^+\rangle.$$
Our goal is to show $\rho (f)\in \TTT_* (l)$ and is a linear
combination of pairs.  To this end,
we consider the following linear combination of toric
quasimodular forms:
$$
\rho_1(f)=\sum_{r+s=k-2}\sum_{m,n=0}^{l-1} \frac {(r+s)!}{r!s!}
\tilde s_{m/l}^{(r+1)}\tilde s_{n/l}^{(s+1)}
\langle f,
(x+y)^ry^s(m,m+n)-(x-y)^ry^s(m,m-n)
\rangle
$$
$$
= C + \sum_{r=0}^{k-2}\sum_{m=0}^{l-1} c_{r,m}\tilde s_{m/l}^{(r+1)} +
\sum_{D>0}q^D\sum_{m,n,r,s} A \langle f,
(x+y)^ry^s(m,m+n)-(x-y)^ry^s(m,m-n) \rangle.
$$
Here $C,c_{r,m}$ are constants whose exact values will not
be needed, and the constant $A=A(r,s,m,D)$ is defined by
$$
A:=\frac {(r+s)!}{r!s!}
\sum_{I (D)}
k_1^rk_2^s
(\delta_{k_1}^{m\mod l}+(-1)^{r-1}\delta_{k_1}^{-m\mod l})
(\delta_{k_1}^{m\mod l}+(-1)^{r-1}\delta_{k_1}^{-m\mod l}),
$$
where $I (D)\subset \ZZ^{4}$ denotes the set
\begin{equation}\label{i.set.def}
I (D) = \{ (m_{1}, k_{1}, m_{2}, k_{2}) \mid  m_1,k_1,m_2,k_2>0, \quad
m_1k_1+m_2k_2=D \}
\end{equation}
A linear combination similar to $\rho_{1} (f)$ appears in the proof of
\cite[Theorem 4.8]{vanish} as the composition of several maps, one of
which is induced by the intersection pairing on Manin symbols.  Here,
however, we just take this as a definition.  After some
simplification, the formula for $\rho_1(f)$ becomes
$$
\rho_1(f)
= C + \sum_{r=0}^{k-2}\sum_{m=0}^{l-1} c_{r,m}\tilde s_{m/l}^{(r+1)}
+ 4\sum_{D>0} q^D
\langle f,
\sum_{I (D)}
(R_{(k_1,k_1+k_2)}^+ - R_{(k_1,k_1-k_2)}^+)
\rangle,
$$
where $R_{(m,n)}^{+}$ is the Manin symbol from Section \ref{rkls}.
The quasimodular form $\rho_1(f)$ and the modular form $\rho(f)$ are
related as follows:

\begin{proposition}\label{firstapprox}
We have
$$
\rho_1(f)-12\rho(f)=
C + 4F_{1} + 4F_{2}+ \sum_{r=0}^{k-2}\sum_{m=0}^{l-1} c_{r,m}\tilde s_{m/l}^{(r+1)},
$$
where $C$ is a constant and 
\begin{align*}
F_{1} &= \sum_{n>0}q^n 
\sum_{d|n}\frac {2n}{d}\langle f,R_{(d,0)}^+\rangle,\\
F_{2} &= \sum_{n>0}q^n 
\sum_{d|n}\Bigl(\langle f,R_{(d,0)}^+\rangle + 2\sum_{\substack{0<e<d}}
\langle f, R_{(d,d-e)}^+ \rangle\Bigr).
\end{align*}
\end{proposition}

\begin{proof}
This follows from the identity
\begin{multline*}
\sum_{I (D)}
(R_{(k_1,k_1-k_2)}^+ - R_{(k_1,k_1+k_2)}^+)
\\ =
-\sum_{d|n}\bigl(\frac {2n}{d}+1\bigr)R_{(d,0)}^+
-2\sum_{\substack{d|n\\d>e>0}}
R_{(d,d-e)}^+
-3\sum_{H (n)}
R_{(c,d)}^+.
\end{multline*}
This identity with weight $k=2$ appears as an intermediate step of the proof of
\cite[Theorem 4.8]{vanish}. However, its proof only uses relations
among $R_{(m,n)}^+$ that are independent of the weight $k$.
\end{proof}

\begin{lemma}\label{f1lemma}
The quasimodular form $F_{1}$ is a linear combination of pairs.
\end{lemma}

\begin{proof}
After some simplification, one can write
\[
F_{1} = 2\sum_{n>0}q^n \sum_{d|n} n d^{k-3} \langle f, x^{k-2}(d,0)_+\rangle.
\]
Let $G$ be the $q$-series 
\[
G  = \sum_{n>0}q^n \sum_{d|n} d^{k-3} \langle f, x^{k-2}(d,0)_+\rangle.
\]
The complex number $\langle f,
x^{k-2}(d,0)_+\rangle$ depends only on $d \mod l$, and further
satisfies 
\[
\langle f, x^{k-2}(-d,0)_+\rangle = (-1)^{k}\langle f, x^{k-2}(d,0)_+\rangle.
\]
Hence $d^{k-3} \langle f, x^{k-2}(d,0)_+\rangle$ is an odd $(\mod
l)$-polynomial.  By Proposition \ref{easyremark}, $G$ is a linear
combination of the $\tilde s_{a/l}^{(k)}$ and a constant, and is hence
toric quasimodular.  Differentiating the linear combination for $G$
with respect to $\tau$ and applying Proposition \ref{derivs.are.there}
completes the proof.
\end{proof}

\begin{lemma}\label{f2lemma}
The quasimodular form $F_{2}+C$ is a linear combination of pairs for a suitably
chosen constant $C$.
\end{lemma}

\begin{proof}
By Proposition \ref{oddprop}, we know that the function 
\[
d\longmapsto\langle f,R_{(d,0)}^+\rangle + 2\sum_{\substack{0<e<d}}
\langle f, R_{(d,d-e)}^+ \rangle
\]
extends to a unique odd $(\mod l)$-polynomial.  The result then
follows from Proposition \ref{easyremark}, and
weight considerations.
\end{proof}

\subsection{}
We are now ready to prove our main theorem.
\begin{theorem}\label{main}
All cusp forms of weight three or more are toric.  Moreover, any such
cusp form can be written as a linear combination of pairs 
(Definition \ref{pairs}).
\end{theorem}

\begin{proof}
One can easily see that lifts of the forms $\tilde s_{a/l}^{(r)}$
can be written as linear combinations of $\,\tilde s\,$ for the new
level.  Therefore, we may assume without loss of generality that $f$
is a newform. Hence by the proof of Proposition \ref{nothingvanishes}, 
$\rho(f)$ is a non-zero multiple of $f$.  Proposition
\ref{firstapprox} and Lemmas \ref{f1lemma} and \ref{f2lemma} show that
$\rho (f)$ can be written up to a constant as a linear
combination of toric quasimodular forms $\tilde s_{a/l}^{(m)}
\tilde s_{b/l}^{(n)}$ for $m+n=k$ and $\tilde s_{a/l}^{(n)}$ for
smaller $n\leq k$. The transformation properties under $\Gamma_1(l)$
insure that all lower weight forms come with zero coefficients, and
that all the quasimodular forms used are actually modular, i.e. 
$E_2s_{a/l}^{(k-2)}$ come with zero coefficients.
\end{proof}

\begin{corollary}
If $l\geq 5$, then any cusp form of weight $k\geq 3$ can be written, up to
a weight $k$ Eisenstein series, as a degree $k$ homogeneous polynomial
in weight one Eisenstein series.
\end{corollary}

\begin{proof}
This follows from Theorem \ref{main}, Proposition \ref{unnamed},
and \cite[Theorem 4.11]{toric}.
\end{proof}

\begin{corollary}
The multiplication map
$$
\MMM_{m}(l)\otimes \MMM_{n}(l) \longrightarrow \MMM_{m+n}(l)
$$
is surjective for all $m\geq n\geq 1$, except for  $m=n=1$.
\end{corollary}

\begin{proof}
Theorem \ref{main} assures that the image of the above map contains
all cusp forms, so it is enough to insure that the forms of the image
take arbitrary values at the cusps. To obtain a form which vanishes
at all but one cusp $p$ we multiply a form in $\MMM_{m}(l)$
that vanishes at all cusps except $p$ and perhaps one other cusp $q$
(relevant only if $m=2$) by a form in $\MMM_{n}(l)$ that
vanishes at $q$ but not at $p$.
\end{proof}

\begin{remark}
A slightly weaker statement can be proved directly by using the fact
that the ring of modular forms is Cohen-Macaulay.  However, we are not
aware of any other proofs for $(m,n)=(2,1)$ or $(2,2)$.
\end{remark}

\begin{remark}
One can also ask which Eisenstein series are toric. It is easy to see
that for a prime level $p$ all Eisenstein series are toric.  For
composite levels, the situation is different.  For example if $l=25$
then weight two toric Eisenstein series form a subspace of codimension
one in the space of all weight two Eisenstein series.  We do not know
any similar examples for higher weight.
\end{remark}

\begin{theorem}\label{eventually}
For every level $l$ there exists an $N$ such that the ring of toric
forms coincides with the ring of modular forms for weights $k\geq N$.
When $l$ is prime, one can take $N=3$.
\end{theorem}

\begin{proof}
In view of Theorem \ref{main}, one needs to show that all Eisenstein
series are eventually contained in the ring of toric forms.  Because
the ring of toric forms is Hecke stable \cite[Theorem 5.3]{toric}, it suffices
to show that the values of toric forms at the cusps eventually
span a $c$-dimensional space, where $c$ is the number of cusps. For
this one needs to show that the values of $s_{a/l}$ for two different
cusps are not proportional.  This is accomplished by a direct
calculation that we leave to the reader.
\end{proof}

\begin{remark}
Theorem \ref{eventually}
was used in \cite{modular} to analyze the embedding of the modular curve $X_{1}
(p)$ given by the graded ring $\TTT_{*} (p)$.
\end{remark}

\section{The map from symbols to forms in higher weight}
\label{s.toricmap}

\subsection{}
A key step in the proof of \cite[Theorem 4.11]{vanish}
was the analysis of a map $\mu$ from the minus space $M_{-}$ of weight 2 Manin
symbols to a quotient of the space $\MMM_{2} (l)$ of weight 2 modular
forms.  Namely, we showed that the map 
\[
\mu \colon (m,n)\longmapsto \tilde s_{m/l} \tilde s_{n/l} 
\]
took $M_{-}$ into the quotient $\MMM_{2} (l)/\EEE_{2} (l)$, where
$\EEE_{2} (l)$ is the space of weight 2 Eisenstein series
(\footnote{This is slightly inaccurate: the map we're denoting
by $\mu$ here is actually the composition of map called $\mu$ in
\cite{vanish} and the Fricke involution.}).  In this
section we consider the analogous map in higher weight given by 
\begin{equation}\label{defofmu}
\mu \colon x^ry^s(m,n)\longmapsto (-1)^s\tilde s^{(s+1)}_{m/l}\tilde s_{n/l}^{(r+1)}
\end{equation}
and describe the relevant quotient containing the image.

\begin{theorem}\label{mumap} 
Let $k>2$.  The map $\mu$ in \eqref{defofmu} applied to the space
generated by the Manin symbols $x^sy^r(m,n)$ takes the relations
\begin{equation}\label{relation}
x^ry^s(a,b)+(-1)^ry^r(x-y)^s(b,-a-b)+(-1)^s(y-x)^rx^s(-a-b,a)
\end{equation}
to the subspace generated by the modular forms $\tilde s_{a/l}^{(k)}$
and the quasimodular forms $\partial_{\tau}\tilde s_{a/l}^{(k-2)}$.
\end{theorem}

\begin{proof}
The symbol \eqref{relation} maps to
\begin{equation}\label{relation2}
(-1)^s\tilde s_{a/l}^{(s+1)}\tilde s_{b/l}^{(r+1)}
+\sum_{t=0}^s\frac{s!}{t!(s-t)!}\tilde s_{b/l}^{(r+t+1)}
\tilde s_{-(a+b)/l}^{(s-t+1)}
+\sum_{t=0}^r\frac{r!}{t!(r-t)!}\tilde s_{-(a+b)/l}^{(r-t+1)}
\tilde s_{a/l}^{(s+t+1)}(-1)^{s+r}.
\end{equation}
Up to quasimodular forms of lower weight
and $\tilde s_{a/l}^{(k)}$, the expression \eqref{relation2} can be
simplified to 
$$
\sum_{D>0}q^D\sum_{I (D)}
(
A_{k_1,k_2}-A_{-k_1,k_2}+A_{k_2,-k_1-k_2} -A_{k_2,k_1-k_2}
+A_{-k_1-k_2,k_1} - A_{k_1-k_2,-k_1}
).
$$
Here $I (D)$ is defined in \eqref{i.set.def} and
$$
A_{k_1,k_2} = (-1)^{s} k_1^sk_2^r\bar\delta^{(a,b)}_{(k_1,k_2)},
$$
where $\bar\delta^{(a,b)}_{(k_1,k_2)}=
\delta_{k_1}^{a\mod l}\delta_{k_2}^{b\mod l}+
(-1)^k\delta_{k_1}^{-a\mod l}\delta_{k_2}^{-b\mod l}$.

The set $I (D)$ can be partitioned into subsets corresponding to 
different ``runs'' of the Euclidean algorithm.  Namely,
there are partially defined maps $\Upsilon $ and $\Delta
$ from $I (D)$ to itself
given by
$$
\Upsilon \colon (m_1,k_1,m_2,k_2)\longmapsto
\left\{\begin{array}{ll}
(m_2,k_1+k_2,m_1-m_2,k_1),&{\rm if~} m_1>m_2\\
(m_2-m_1,k_2,m_1,k_1+k_2),&{\rm if~} m_1<m_2\\
{\rm not~defined},&{\rm if~} m_1=m_2
\end{array}\right.
$$
$$
\Delta \colon (m_1,k_1,m_2,k_2)\longmapsto
\left\{\begin{array}{ll}
(m_1+m_2,k_2,m_1,k_1-k_2),&{\rm if~} k_1>k_2\\
(m_2,k_2-k_1,m_1+m_2,k_1),&{\rm if~} k_1<k_2\\
{\rm not~defined},&{\rm if~} k_1=k_2
\end{array}\right.
$$
These maps are inverses of each other whenever their composition is
defined.  The whole set $I (D)$ can be pictured as a disjoint union of
vertical threads, where each thread is obtained by starting at the top
with a solution with $m_1=m_2$ and applying $\Delta $
until arriving at a solution with $k_1=k_2$(\footnote{$\Upsilon$ and
$\Delta$ stand for \emph{up} and \emph{down}.}).  The crucial
observation is that for each thread $\Theta $, the sum 
$$
\sum_{\Theta } A_{k_1,k_2} + A_{k_2,-k_1-k_2} + A_{-k_1-k_2,k_1} - A_{-k_1,k_2} - A_{k_2,k_1-k_2}
-A_{k_1-k_2,-k_1}
$$
collapses.  Indeed, the negative terms
for elements $(m_1,k_1,m_2,k_2)$ cancel the positive terms for
elements $\Delta (m_1,k_1,m_2,k_2)$.  To see this, observe that if
$k_1>k_2$, then the positive terms of $\Delta (m_1,k_1,m_2,k_2)$ equal
$$
A_{k_2,k_1-k_2}+A_{k_1-k_2,-k_1}+A_{-k_1,k_2}.
$$
The $k_1<k_2$ case is handled similarly, taking into account the symmetry
$A_{-k_1,-k_2}=A_{k_1,k_2}$.

Hence, up to a
linear combination of lower weight forms and the forms $\tilde
s_{a/l}^{(k)}$, the image of the relation \eqref{relation} is equal to
\begin{multline*}
\sum_{D>0}q^D
\Bigr(\sum_{\{i\in I (D) \mid m_1=m_2\}}(A_{k_1,k_2}+A_{k_2,-k_1-k_2}+A_{-k_1-k_2,k_1})
\\-\sum_{\{i\in I (D) \mid k_1=k_2\}}(A_{-k_1,k_2}+A_{k_2,k_1-k_2}+A_{k_1-k_2,-k_1})
\Bigl).
\end{multline*}
The coefficient of $q^{D}$ can be further simplified to 
\begin{multline*}
\sum_{d|D}\sum_{0<e<d}
(A_{e,d-e}+A_{d-e,-d}+A_{-d,e})
-\sum_{d|D}(\frac Dd -1) \Bigl(
d^{k-2}(\delta_{(-d,d)}^{(a,b)} +(-1)^k\delta_{(-d,d)}^{(-a,-b)})
\\+(-1)^sd^s0^r(\delta_{(d,0)}^{(a,b)} +(-1)^k\delta_{(d,0)}^{(-a,-b)})
+0^sd^r(\delta_{(0,d)}^{(a,b)} +(-1)^k\delta_{(0,d)}^{(-a,-b)})
\Bigr),
\end{multline*}
where $\delta$ is now a Kronecker symbol for elements of $(\ZZ/l\ZZ)^2$,
and our convention is $0^s=1$ if and only if $s=0$.

To finish the proof we first observe that the contribution of the
terms with $D/k$ is, up to an additive constant,
a derivative with respect to $\tau$ of
$$
-\delta_{a+b}^0\tilde s_{b/l}^{(k-2)}+(-1)^{s+1}0^r\delta_{b}^0
\tilde s_{a/l}^{(k-2)}
-0^s\delta_{a}^0\tilde s_{b/l}^{(k-2)},
$$
where $\delta$ is the usual Kronecker function.
To show that the remaining contributions give linear combinations of
the forms $\tilde s_{a/l}^{(\leq k)}$, it is enough to establish that
for any $a,b,r,s$ the $(\mod l)$-polynomial
$$h(d):=\sum_{0<e<d}
(A_{e,d-e}+A_{d-e,-d}+A_{-d,e})
+ A_{-d,d}+A_{d,0}+A_{0,d}
$$
is odd. This follows easily from Proposition \ref{conesum} and the
symmetry of $A$.
\end{proof}

\begin{corollary}\label{truemu}
The map $\mu$ induces a map from the space of weight $k$ Manin symbols
$M$ to the quotient $\QQQ$ of the space of weight $k$ quasimodular forms by
subspace generated by the Eisenstein series $\tilde s_{a/l}^{(k)}$ and
the derivatives $\partial_\tau \tilde
s_{a/l}^{(k-2)}$.
\end{corollary}

\begin{remark}\label{toric1}
An alternative approach to Theorem \ref{mumap} is to
look at the identify
$$
(s_\alpha^{(1)}+s_\beta^{(1)}+s_{-\alpha-\beta}^{(1)})^2
+\frac 12(s_\alpha^{(2)}+s_\beta^{(2)}+s_{-\alpha-\beta}^{(2)})=0
$$
which comes from a calculation of certain toric form for
the complex projective plane $\PP^2$, see \cite{toric}.
One can differentiate the above identity with respect to
$\alpha$ and  $\beta$ several times and  plug in rational values of
$\alpha$ and $\beta$. Then it remains to use the transformation
that connects $\tilde s_{a/l}^{(k)}$ and $s_{a/l}^{(k)}$. We
leave the details to the reader.
\end{remark}

\section{Hecke equivariance of the symbols to forms map}\label{s.Hecke}

\subsection{}
It is not hard to see by explicit computation that the subspace
spanned by the Eisenstein series and derivatives mentioned in
Corollary \ref{truemu} is invariant under the action of $\Gamma_0 (l)
/\Gamma_1 (l)$, the Fricke involution, and the Hecke operators.  Hence
we can naturally extend their action to the quotient $\QQQ$.  The goal
of this section is to show that the map of Corollary \ref{truemu} is
compatible with the action of Hecke operators.  For this, one needs to
show that the map
$$
x^ry^s(m,n)\mapsto (-1)^s\tilde s^{(s+1)}_{m/l}\tilde s_{n/l}^{(r+1)}
$$
is compatible with the action of Hecke operators, up to linear combinations
of $\tilde s_{a/l}^{(k)}$ and $\partial_\tau \tilde s_{a/l}^{(k-2)}$.

\begin{theorem}\label{HE}
Let $p$ be a prime number coprime to $l$ and $T_p$ be the
corresponding Hecke operator on $M_k$ and $\MMM_k$, where we abuse
notations slightly.  Let $\mu$ be the map defined in Theorem
\ref{mumap}.  Then for every $w\in M_k$, the image $\mu(T_p w)$ is
equal to $T_p(\mu(\epsilon_{p^{-1}}w))$ modulo a linear combination of
$s_{a/l}^{(k)}$ and $\partial_\tau s_{a/l}^{(k-2)}$.  Here
$\epsilon_{p^{-1}}$ is the action of the element of
$\Gamma_0(l)/\Gamma_1(l)$ given by $x^ry^s(u,v)\mapsto x^ry^s(pu,pv)$,
(cf. Proposition \ref{gamma0action}).
\end{theorem}

Before we begin the proof of Theorem \ref{HE}, we need a lemma
giving a geometric interpretation of the set $H (p)$ involved 
in Merel's description of the $T_{p}$-action on Manin symbols
(Theorem \ref{HeckeisHecke}).

\begin{lemma}\label{abcd}
\cite[Theorem 3.16]{vanish} For each index $p$ sublattice $S\subset
\ZZ^{2}$, consider the convex hull of all nonzero points of $S$ that
lie in the first quadrant. Then the compact subset of the boundary of
this convex hull is a union of segments.  Moreover, the coordinates
$(a,c)$, $(b,d)$ of the vertices of each segment (ordered from the
$x$-axis) satisfy $ad-bc=p$ and $a>b\geq 0$, $d>c\geq 0$, and hence
determine an element of $H (p)$.  Conversely, all $(a,b,c,d)\in H (p)$
come from one such sublattice $S$ in this manner.
\end{lemma}

Given an index $p$ sublattice $S\subset \ZZ^{2}$, we write $H
(p,S)$ for the subset of those $(a,b,c,d)\in H (p)$ corresponding to
$S$.

\begin{example}
Figure \ref{Fig1} shows the case $p=2$.  There are three
sublattices of index $2$, and altogether four distinct boundary segments.
>From the segments we obtain the four elements of $H (2)$,
namely $(1,0,0,2)$, $(2,1,0,1)$, $(1,0,1,2)$ and $(2,0,0,1)$.
\begin{figure}[tbh]
\begin{center}
\includegraphics[scale = 0.4]{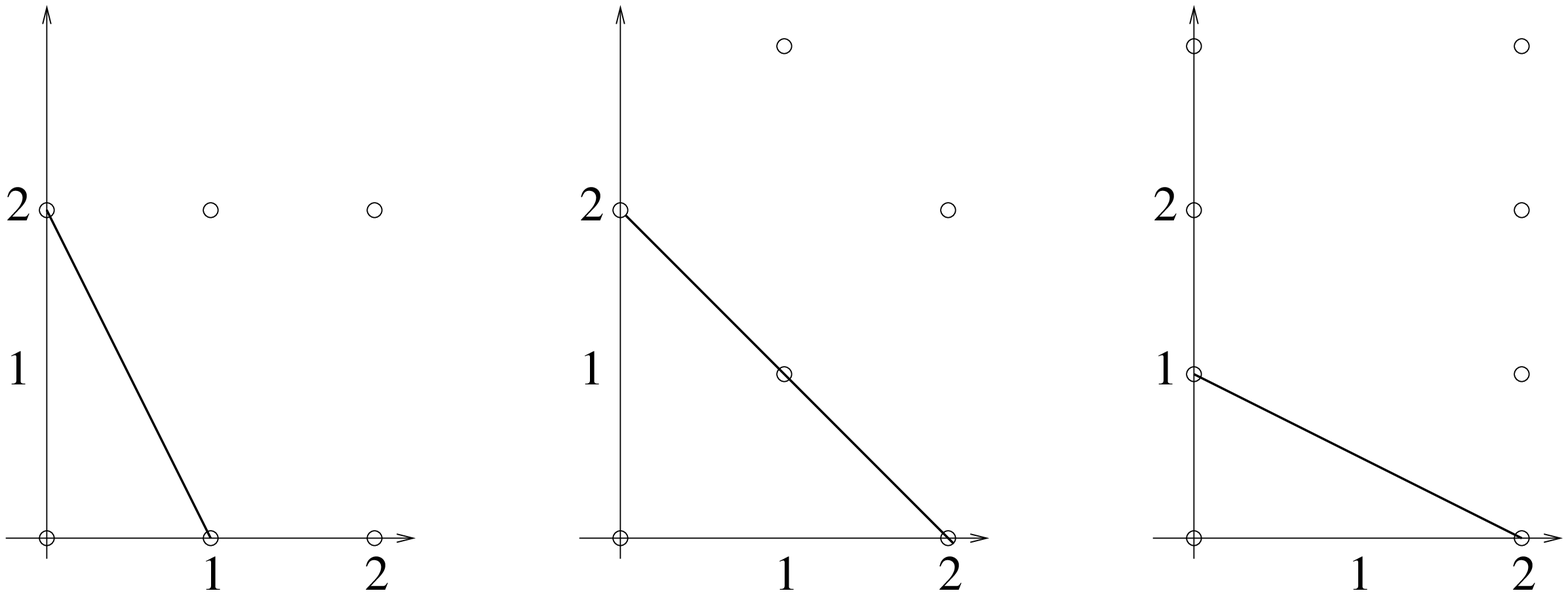}
\end{center}
\caption{\label{Fig1}}
\end{figure}
\end{example}

We will also need the following duality operation on the set of sublattices.

\begin{definition}
For an index $p$  sublattice $S$ we denote by $S^*$ the
sublattice of all points $P$ in $\ZZ^2$ such that $P\cdot S\subseteq p\ZZ$.
where $\cdot$ is the standard scalar product on $\ZZ^2$.
It is clear that $S^{**}=S$. Moreover, $S^*$ can be obtained from $S$
by a $\pi/2$ rotation at the origin.
\end{definition}

We are now ready to start the proof of Theorem \ref{HE}.

\begin{proof}
It is enough to consider $w= x^ry^s(u,v)$.
By Theorem \ref{HeckeisHecke} and the definition of $\mu$,
\[
\mu(T_p x^ry^s(u,v))\sim_l
\sum_Dq^D\sum_{h\in H (p)}
\sum_{i\in I (D)} \Phi (h,i), 
\]
where 
\[
\Phi (h,i) = 
(ak_2-bk_1)^r(ck_2-dk_1)^s
\bar\delta_{k_1,k_2}^{au+cv,bu+dv} -
(ak_2+bk_1)^r(ck_2+dk_1)^s
\bar\delta_{k_1,k_2}^{-au-cv,bu+dv}.
\]
Here $\sim_l$ means that equality holds modulo linear combinations
of $\tilde s_{a/l}^{(<k)}$, and $I (D)$ is defined in
\eqref{i.set.def}.
We can use $(p,l)=1$ to rewrite the above
as
$$
\Phi (h,i) =  A_{dk_1-ck_2,ak_2-bk_1}-A_{-dk_1-ck_2,ak_2+bk_1},
$$
where $A_{\alpha,\beta}=\beta^r(-\alpha)^s\bar\delta_{\alpha,\beta}^{pu,pv}$.
On the other hand,
\begin{multline}\label{first.stab}
T_p \mu(\epsilon_{p^{-1}}w)\sim_{pl}
\sum_Dq^D\sum_{I (pD)}
(A_{k_1,k_2}-A_{-k_1,k_2})
\\
+p^{k-1}
\sum_Dq^D\sum_{I (D)}
(-1)^sk_1^sk_2^r(\bar\delta_{k_1,k_2}^{u,v}-\bar\delta_{-k_1,k_2}^{u,v}).
\end{multline}
For each $i\in I (pD)$ there exists a sublattice $S$ such that $(m_1,m_2)\in S$ and
$(k_1,k_2)\in S^*$. Moreover, $S$ is unique unless $m_1,k_1,m_2,k_2
=0\mod p$, in which case there are $(p+1)$ such sublattices $S$.
To record this, we use the notation 
\[
I (pD,S) = \{i\in I (pD) \mid  (m_1,m_2)\in S, \quad (k_{1},k_{2})\in S^{*}\}.
\]
Let us further write, for any two subsets $U_{1}, U_{2}\subset \RR^{2}$,
\[
I(pD,S;U_{1}, U_{2}) = \{i\in I (pD,S) \mid  (m_1,m_2)\in U_{1}, (k_{1}, k_{2})\in U_{2}\}.
\]
Now we can rewrite \eqref{first.stab} as
\begin{equation}\label{one.guy}
T_p \mu(\epsilon_{p^{-1}}w)\sim_{pl}
\sum_Dq^D
\sum_S
\Bigl(
\sum_{I (pD,S;Q_{I}, Q_{I})} A_{k_1,k_2}
-
\sum_{I (pD,S;Q_{II}, Q_{II})} A_{k_1,k_2}
\Bigr)
\end{equation}
where $Q_{I}$ and $Q_{II}$ denote the open first and the second quadrants.
\begin{remark}
The reason we must write $\sim_{pl}$ here rather than $\sim_l$ is that
the action of $T_p$ defined for weight $k$ on $\tilde s_{a/l}^{(\leq k)}$
will be a linear combination $\tilde s_{a/pl}^{(\leq k)}$.
\end{remark}

Given any $h = (a,b,c,d)\in H (p)$, we also denote by $h$ the linear
transformation $\RR^{2}\rightarrow \RR^{2}$ given by the multiplying by 
matrix $\left
( \begin{smallmatrix}a&b\\c&d\end{smallmatrix}\right)$ on the right.  
This allows
us to write $\mu(T_p(w))-T_p \mu(\epsilon_{p^{-1}}w)$ as
\begin{equation}\label{second.stab}
\mu(T_p(w))-T_p \mu(\epsilon_{p^{-1}}w)\sim_{pl}
\sum_{D}q^D\sum_S
(\Sum_{1} - \Sum_{2} - \Sum_{3} + \Sum_{4}),
\end{equation}
where 
\begin{align*}
\Sum_{1} &= \sum_{\substack{h\in H (p,S)\\
i\in I (pD,S;h^t (Q_{I}), h^{-1} (Q_{I}))}} A_{k_1,k_2}\\
\Sum_{2} &= \sum_{\substack{h\in H (p,S)\\
i\in I (pD,S;h^t (Q_{II}), h^{-1} (Q_{II}))}}A_{k_1,k_2}\\
\Sum_{3} &= \sum_{I (pD,S,Q_{I})} A_{k_1,k_2}\\
\Sum_{4} &= \sum_{I (pD,S,Q_{II})} A_{k_1,k_2}
\end{align*}
It is convenient to visualize $\Sum_{1},\dotsc ,\Sum_{4}$ as indicated
in Figure \ref{Fig2}.

\begin{figure}[tbh]
\begin{center}
\includegraphics[scale = 0.4]{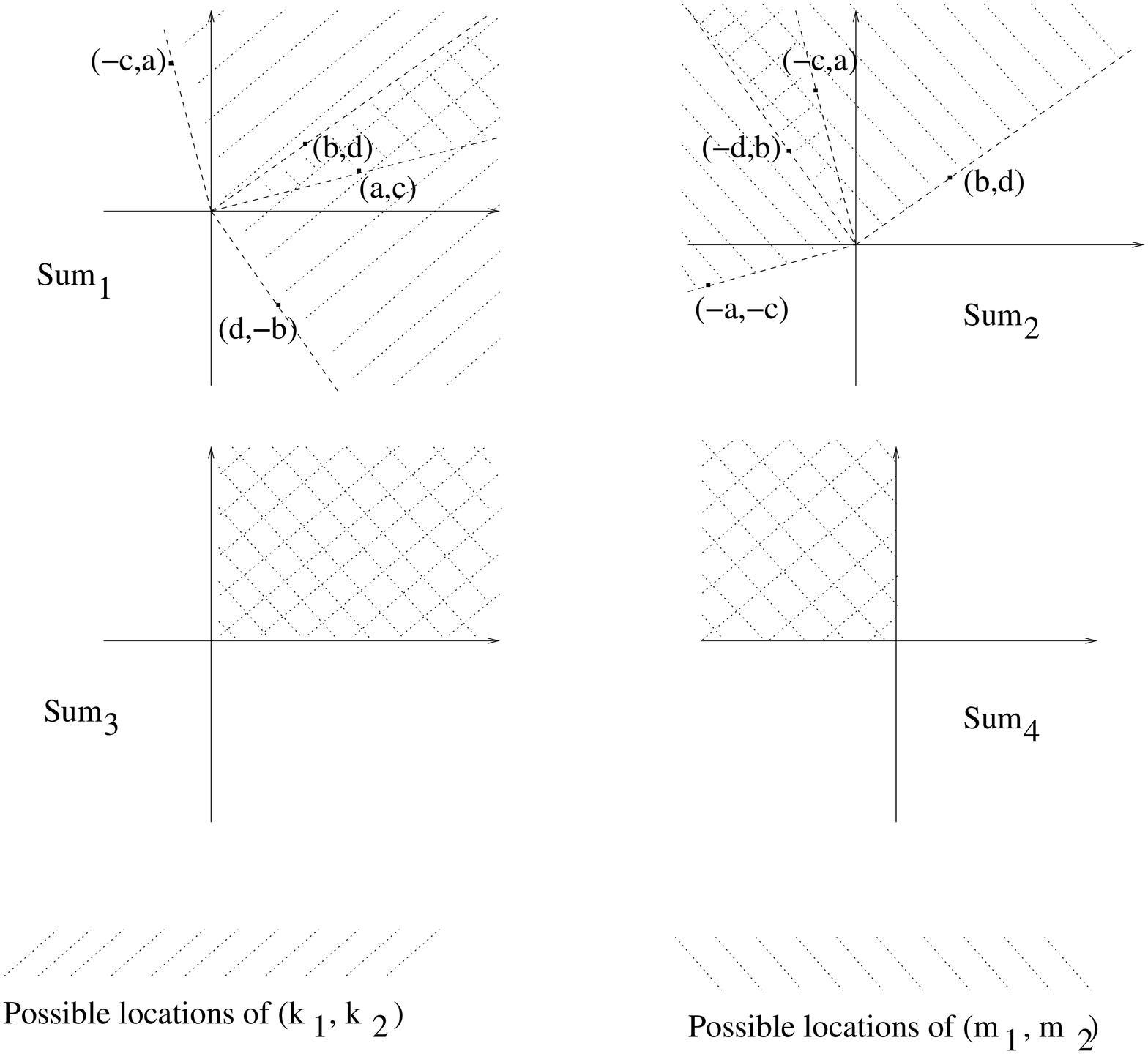}
\end{center}
\caption{\label{Fig2}}
\end{figure}

\subsection{}
Now we consider the right hand side of \eqref{second.stab}.
It will turn out that most terms will cancel each other, but there
will be some terms left over that will require careful consideration.
To discuss these terms, we require some additional terminology.

For every sublattice $S$ of index $p$ we form the convex hulls of the
nonzero points in \emph{each} quadrant.  The open 1-cones spanned by
the points on the boundary of these hulls will be called {\em rays} of
the lattice $S$, and will be denoted $\rho (S)$.  The rays generated
by the points $(\pm p,0)$ and $(0,\pm p)$ will be called the
\emph{axis} rays; all others will be called \emph{non-axis} rays.  By
abuse of notation, given a point $(x,y)$ we write $(x,y)\in \rho (S)$
to mean that $(x,y)$ lies on a ray of $S$.  
Finally, 
for any nonzero point $v\in S$, we define a rational cone $\cone(v)$
as follows:
\begin{itemize}
\item If $v$ lies on a ray $\rho$, then we put $\cone (v) = \rho$.
\item Otherwise, we set $\cone (v)$ to be the interior of the unique $2$-cone 
spanned by adjacent rays of $S$ and containing $v$.
\end{itemize}

The first step in investigating \eqref{second.stab} is
the following lemma, which is the heart of the proof.
\begin{lemma}\label{inside}
Let $S\subset \ZZ^{2}$ have index $p$.  Then for every
$(m_1,k_1,m_2,k_2)$ such that $(m_1,m_2)\not \in \rho (S)$ and
$(k_1,k_2)\not \in \rho (S^{*})$, the total of the contributions
of $\Sum_1$, $\Sum_2$, $\Sum_3$ and $\Sum_4$ at $(m_1,k_1,m_2,k_2)$
and $-(m_1,k_1,m_2,k_2)$ is zero.
\end{lemma}

\begin{proof}[Proof of Lemma \ref{inside}]
Clearly, it is enough to check this if $(k_1,k_2)$ is in the first or
second quadrant.

First assume $(k_1,k_2)\in Q_{I}$.  Then the only nontrivial
contributions come from $\Sum_1$ and $\Sum_3$ when $(m_1,m_2)\in
Q_{I}$, and in this case we claim $\Sum_{1}$ contributes $A_{k_1,k_2}$
and $\Sum_3$ contributes $-A_{k_1,k_2}$.  Indeed, if $(m_1,m_2)\in
Q_{I}$, there is a contribution of exactly one $(a,b,c,d)$ in
$\Sum_1$, which corresponds to $\cone(m_1,m_2)$.  Hence the total
contribution is zero.

Next assume $(k_1,k_2)\in Q_{II}$.  In this case $\Sum_{3}$ doesn't
contribute, and we split the contributions of the remaining sums into types:
\begin{itemize}
\item ($\Sum_1$, type $1$)  We assume $(m_1,m_2)\in Q_{I}$
lies in a cone above $\cone (k_2,-k_1)$
(see Figure \ref{Fig3}, graph 1). Then there is a unique $(a,b,c,d)$ in $\Sum_1$,
that corresponds to $\cone(m_1,m_2)$, and the contribution of $\Sum_1$
is $A_{k_1,k_2}$.
\item ($\Sum_1$, type $2$) We assume $(m_1,m_2)\in Q_{I}$ lies
in a cone below the $\cone (k_2,-k_1)$
(see Figure \ref{Fig3}, graph 2).
Again, there is one $(a,b,c,d)$, and the contribution is
$A_{-k_1,-k_2}$.
\item ($\Sum_2$, type $1$)  We assume $(m_1,m_2)\in Q_I$ 
lies in a cone above $\cone (k_2,-k_1)$.
Then there is a unique $(a,b,c,d)$ in $\Sum_2$ that corresponds to $\cone(k_2,-k_1)$,
and the contribution of $\Sum_2$ is $-A_{k_1,k_2}$. 
\item ($\Sum_2$, type $2$)  We assume $(m_1,m_2)\in Q_{I}$ such that $(-m_1,-m_2)$
lies in a cone below $\cone (k_2,-k_1)$.
Then there is a unique $(a,b,c,d)$ in $\Sum_2$ that corresponds to $\cone(k_2,-k_1)$,
and the contribution of $\Sum_2$ is $-A_{k_1,k_2}$.
\item ($\Sum_2$, type $3$) If $(m_1,m_2)\in Q_{II}$, then there is a contribution
of $-A_{k_1,k_2}$. Indeed, the unique $(a,b,c,d)$ corresponds to
$\cone(k_2,-k_1)$.
\end{itemize}
Clearly, the type $1$ contributions of $\Sum_{1}$ and $\Sum_2$ cancel;
after we apply the symmetry $A_{-k_1,-k_2}=A_{k_1,k_2}$, the type $2$
contributions of $\Sum_{1}$ and $\Sum_2$ cancel as well.  Finally, the
contribution of $\Sum_4$ cancels the type $3$ contribution of
$\Sum_2$, which completes the proof of Lemma \ref{inside}.
\end{proof}

\begin{figure}[tbh]
\begin{center}
\includegraphics[scale = 0.4]{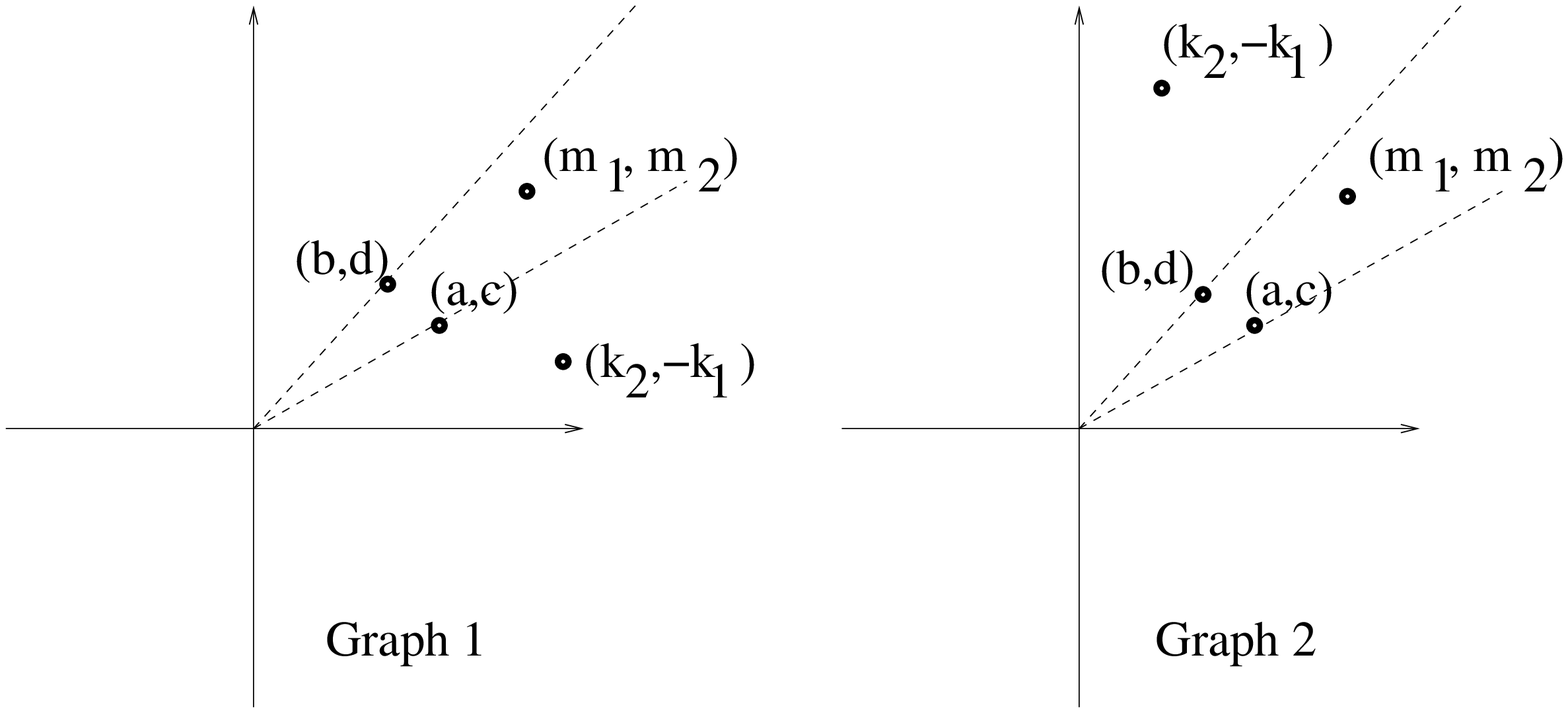}
\end{center}
\caption{\label{Fig3}}
\end{figure}

\begin{remark}
In the terms that cancel, the matrices $(a,b,c,d)$ are different,
which makes Lemma \ref{abcd} crucial to the success of the proof.
\end{remark}

We now return to the proof of Theorem \ref{HE}. Having handled the
bulk of the terms in $\Sum_1$--$\Sum_4$, we now examine the
cases when $(m_1,m_2)$ or $(k_1,k_2)$ lie on a ray.  For any point
$(u,v)$, we let $(u,v)^{\perp}$ be the set of all $(x,y)$ with $ux+vy=0$.  

For each $(m_{1},m_{2})\in \rho (S)$ we define a subset $C
(m_{1},m_{2})\subset S^{*}$ as follows.  If $(m_{1},m_{2})$ is not on
a coordinate axis, then we let $C (m_{1},m_{2})$ be the set of all
points with positive scalar product with $(m_{1},m_{2})$ \emph{except}
those that lie in one of the closed cones adjacent to
$(m_{1},m_{2})^{\perp}$ (Figure \ref{Fig4}).  We use the same notation
to denote the similar set $C (k_{1},k_{2})$ constructed from a point
$(k_{1},k_{2})\in \rho (S^{*})$.  If $(m_{1},m_{2})$ or $(k_{1},
k_{2})$ lies on an axis, then we define $C (m_{1},m_{2})$ and $C
(k_{1},k_{2})$ using the small diagrams in Figure \ref{Fig4}.

\begin{figure}[tbh]
\begin{center}
\includegraphics[scale = 0.4]{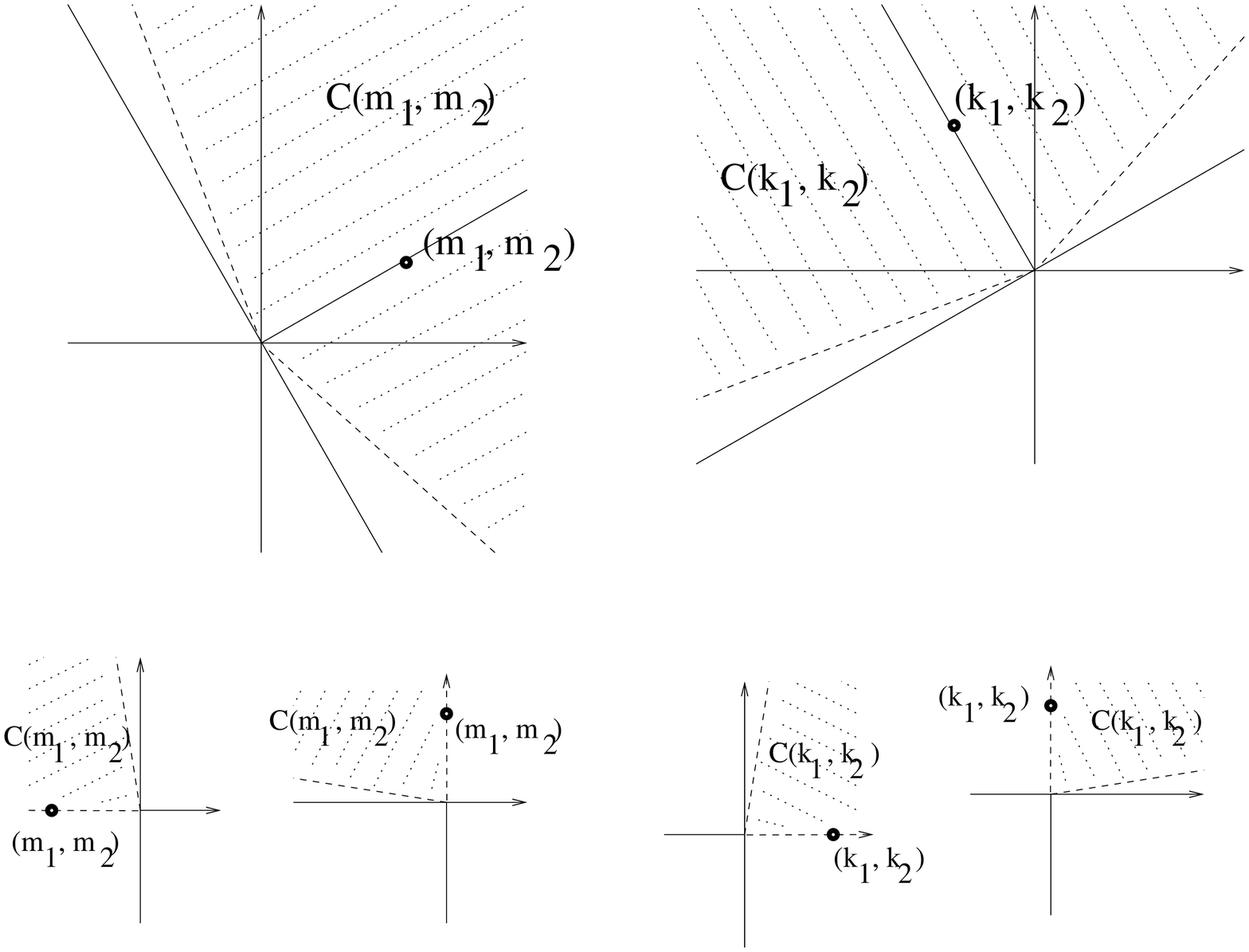}
\end{center}
\caption{\label{Fig4}}
\end{figure}

The rays of the boundary of $C(m_1,m_2)$ (excluding the origin) will
be denoted by $\partial C(m_1,m_2)$ and similarly for $\partial
C(k_1,k_2)$.  For any cone $C$, we write $\sideset{}{'}\sum_C$ to indicate
that the sum
is taken over $C\cup \partial C$ with terms lying in $\partial C$ taken
with weight $1/2$.

\begin{lemma}\label{rayslemma}
With the above notation,
\begin{multline}\label{rayslemmaeq}
\mu(T_p(w))-T_p \mu(\epsilon_{p^{-1}}w)
\sim_{pl}\\
\sum_S
\sum_{(k_1,k_2)\in \rho (S^*)\cap Q'_{II}}
\sideset{}{'}\sum_{(m_1,m_2)\in C(k_1,k_2)}
q^{(m_1k_1+m_2k_2)/p}A_{k_1,k_2}\\
-
\sum_S\sum_{(m_1,m_2)\in \rho (S)\cap Q'_{I}}
\sideset{}{'}\sum_{(k_1,k_2)\in C(m_1,m_2)}
q^{(m_1k_1+m_2k_2)/p}A_{k_1,k_2},
\end{multline}
where $Q'_{I}$ and $Q'_{II}$ are the closures of the first and second
quadrants.
\end{lemma}

\begin{proof}
[Proof of Lemma \ref{rayslemma}]
Because of Lemma \ref{inside},
we need to examine the contribution of $\Sum_1$, $\Sum_2$, $\Sum_3$ and $\Sum_4$
to the quadruples $\pm(m_1,k_1,m_2,k_2)$ where at least one of $(m_1,m_2)$
and $(k_1,k_2)$ lie on the ray of the corresponding lattice.
We will have to be especially careful when one of these vectors is
located on a coordinate axis. In what follows we will fix lattices 
$S$ and $S^*$.

First, let us deal with the case when $(k_1,k_2)\in \rho (S^*)$ and is
not on an axis, and $(m_1,m_2)\not \in \rho (S)$.  If $(k_1,k_2)$ is
in the first or third quadrant, then $\Sum_2$ and $\Sum_4$ do not
contribute, and the contributions of $\Sum_1$ and $\Sum_3$ cancel
since they are respectively $A_{k_1,k_2}$ and $-A_{k_1,k_2}$.
Therefore, it is enough to consider when $(k_1,k_2)$ lies in the
second or fourth quadrants.
Now the terms of $\Sum_2$ and $\Sum_3$ do not
contribute, and the total contribution of $\Sum_1$ and
$\Sum_4$ equals $A_{k_1,k_2}$ if and only if $(m_1,m_2)\in
C(k_1,k_2)$, and is zero otherwise.  This clearly corresponds to the
terms we get on the right of \eqref{rayslemmaeq}.  

Now suppose $(k_1,k_2)$ lies on a coordinate axis.  We may assume that
it lies in the positive portion.  If $(m_1,m_2)\not \in \rho (S)$,
then only $\Sum_1$ contributes, and the contribution is $A_{k_1,k_2}$
if any only if $(m_1,m_2)\in C(k_1,k_2)$. Clearly this corresponds
exactly to the contribution on the right of \eqref{rayslemmaeq}.

Analogously one can treat the case of $(m_1,m_2)\in \rho (S)$
with $(k_1,k_2)\not \in \rho (S^*)$. We have therefore shown that
Lemma \ref{rayslemma} holds up to the contributions of $\pm(m_1,k_1,m_2,k_2)$
with both $(m_1,m_2)$ and $(k_1,k_2)$ on the rays of the corresponding
lattices.

If both $(m_1,m_2)$ and $(k_1,k_2)$ belong to non-axis rays of
$S$ and $S^*$, the contributions of $\Sum_1$ and $\Sum_2$ are zero. Hence,
the contribution of $-A_{k_1,k_2}$ occurs if both of them lie in the
first or third quadrant and the contribution of $A_{k_1,k_2}$ occurs
if both lie in the second or fourth quadrant. To show that this
is consistent with the right hand side of the equation
of the lemma, observe that if $(k_1,k_2)\in Q_{II}$ and $(m_1,m_2)\in Q_I$,
the contributions of the two $\sideset{}{'}\sum$ cancel. Indeed, in this
case $(k_1,k_2)\in C(m_1,m_2)$ and $(k_1,k_2)\in \partial C(m_1,m_2)$
is equivalent to $(m_1,m_2)\in C(k_1,k_2)$ and
$(m_1,m_2)\in \partial C(k_1,k_2)$, respectively.

The remaining case of one or both of $(m_1,m_2)$ and $(k_1,k_2)$
on the axis with both of them on the rays is treated similarly and is left
to the reader.
\end{proof}

Continuing now with the proof of Theorem \ref{HE}, we investigate the
sums on the right of \eqref{rayslemmaeq}.  We divide the contributions
to the sums over $\rho (S^{*})$ into two types: those coming from
non-axis rays, and those coming from axis rays.

\begin{lemma}\label{nonaxislemma}
In the sums over $S$ in \eqref{rayslemmaeq}, the contributions of the
non-axis rays give a linear combination of $\tilde s_{a/pl}^{(\leq
k)}$ and $\partial_\tau \tilde s_{a/l}^{(k-2)}$.
\end{lemma}

\begin{proof}
[Proof of Lemma \ref{nonaxislemma}]
First we calculate the
contribution of a $(k_1,k_2)\in \rho (S^*)$ such that $(k_1,k_2)$ lies
on the ray $\RR_{>0}(-c,a)$, where $(a,c)$ is in the first quadrant.

Let $(b,d)$ (respectively $(b_1,d_1)$) be the generator of the ray of
$S$ adjacent to the ray generated by $(a,c)$ in the counterclockwise
(resp. clockwise) direction.  Then the sets of vectors $\{(a,c), (b,d)
\}$ and $\{(a,c), (b_{1}, d_{1}) \}$ form a $\ZZ$-basis of $S$, 
which implies $(b,d)+(b_1,d_1)=N(a,c)$ where $N$ is a positive integer.
Then any $(m_1,m_2)\in C(k_1,k_2)$ can be written
\[
(m_1,m_2)=-\alpha(a,c)+\beta(b,d),
\]
where
\begin{equation}\label{conditions}
\alpha,\beta\in\ZZ,\quad (m_1,m_2)\cdot(-d,b)>0,\quad (m_1,m_2)\cdot(-d_1,b_1)>0.
\end{equation}
The conditions \eqref{conditions} translate into the inequality
$0<\alpha<N\beta$
on $\alpha$, which has $N\beta-1$ solutions for a given $\beta$.
Note that the terms in $\partial C(K_1,K_2)$ correspond to $\alpha=0$ and
$\alpha= N\beta$, which contributes an extra $A_{k_1,k_2}$
for each value of $\beta$.

Now if we write $(k_1,k_2)=t(-c,a)$ for some positive integer $\ZZ $,
then $(m_1k_1+m_2k_2)/p=t\beta$, so that 
the contribution of the complete ray $\RR_{>0}(-c,a)\in \rho (S^*)$ to
the first term of Lemma \ref{rayslemma} is
$$
\sum_{t>0} \sum_{\beta>0} q^{t\beta} (N\beta-1)A_{-tc,ta}=
\sum_{D>0} q^D\sum_{t|D} \frac {ND}t A_{-tc,ta}.
$$
When one recalls the definition of $A$, this is easily
seen to be a linear combination of
$\partial_\tau\tilde s_{a/l}^{(k-2)}$.

Next we calculate the contribution of an $(m_1,m_2)$ that lies on ray
$\RR_{>0}(a,c)$ of $S$.  The computation is very similar to the
above.  As before we denote by $(b,d)$ and
$(b_1,d_1)$ the generators of the rays of $S$ adjacent to
$\RR_{>0}(a,c)$. Then in the second summation of \eqref{rayslemmaeq},
the pairs $(k_1,k_2)$ are of the form
$$(k_1,k_2)=-\alpha(c,-a)+\beta(d,-b),~\alpha,\beta\in\ZZ,$$
where as before $0<\alpha<N\beta$ for $(m_1,m_2)\in C(m_1,m_2)$
and $\alpha=0$ or $\alpha =N\beta$ for $(m_1,m_2)\in \partial C(m_1,m_2)$.
If we write $(m_1,m_2)=t(a,c)$ for $t$ a positive integer, then we 
obtain
$$
-\sum_{t>0}\sum_{\beta>0}q^{t\beta}
\Big(
\sum_{0<\alpha<N\beta}A_{-\alpha c+\beta d,\alpha a -\beta b}
+\frac12
A_{\beta d, -\beta b}
+\frac12
A_{\beta(-Nd+c), \beta(Na-b)}
\Big).
$$
It remains to use Propositions \ref{conesum} and \ref{easyremark}
to see that the above is a linear combination of $s_{a/l}^{(\leq k)}$.
This completes the proof of the lemma.
\end{proof}

\begin{lemma}\label{axislemma}
In the sums over $S$ in \eqref{rayslemmaeq}, the contributions of the
axis rays give a linear combination of $\tilde s_{a/pl}^{(\leq
k)}$ and $\partial_\tau \tilde s_{a/l}^{(k-2)}$.
\end{lemma}

\begin{proof}
[Proof of Lemma \ref{axislemma}]
First, if $S$ or $S^*$
contain $(0,1)$ or $(1,0)$, then the contributions of the two sums in
Lemma \ref{rayslemma} cancel.  Hence we may ignore lattices of this type.

If $(k_1,k_2)$ is on the positive half
of the $x$-axis, then $k_1$ is a multiple of $p$.
The top cone of $S$ in the first quadrant is
the span of the positive half of $y$-axis and $(1,a)$, with $a$ taking
all values from $0$ to $p-1$, depending on $S$. One then observes that
the contribution of $S_1$ and $S_2$ with $a_1+a_2=p$ can be thought of
as the sum over $(m_1,m_2)$ in the interior of the cone spanned by
$(1,a_1)$ and $(1, -a_2)$, plus half the sum for $(m_1,m_2)$ on
the boundary of the cone. It is then easily seen to give
a linear combination of $\partial_\tau\tilde s_{a/l}^{(k-2)}$.
The case of $(k_1,k_2)$ on the positive half of the $y$-axis
is treated similarly.

If $(m_1,m_2)$ is on one of the axes, then we observe that the sum of
$A_{k_1,k_2}$ over $C(m_1,m_2)$ and its boundary can be thought of as
the sum over all points of $\ZZ^2$ that lie in that cone of an even
two-variable $\mod pl$-polynomial $\hat A_{k_1,k_2}$, which we define
to equal $A_{k_1,k_2}$ if $(k_1,k_2)\in S^*$ and zero otherwise. One
then again invokes Propositions \ref{conesum} and \ref{easyremark} to
conclude that these terms contribute a linear combination of $\tilde
s_{a/pl}^{(\leq k)}$.
\end{proof}

\noindent\textit{Completion of the proof of Theorem \ref{HE}}.
By Lemmas \ref{nonaxislemma} and \ref{axislemma}, we have that 
$$
\mu(T_p(w))-T_p \mu(\epsilon_{p^{-1}}w)
$$
is a linear combination of $\partial_\tau \tilde s_{a/l}^{(k-2)}$
and $\tilde s_{a/pl}^{(\leq k)}$.   The modular transformation properties then
imply that only $\tilde s_{a/l}^{(k)}$ and
$\partial_\tau\tilde s_{a/l}^{(k-2)}$ appear, which
finishes the proof of Theorem \ref{HE}.
\end{proof}

\begin{remark}
Another way to state Theorem \ref{HE} is to say that the composition of
$\mu$ and Fricke involution is Hecke-equivariant.
\end{remark}

\begin{remark}\label{toric2}
The discussion of this section simplifies a bit if one uses the
geometry of toric varieties. More specifically, one has to
consider toric modular forms $f_{\ZZ^2,\deg}$ defined in \cite{toric}
and then differentiate them with respect to the components of the
degree function $\deg$. Then the Hecke action described in \cite{toric}
can be interchanged with these partial differentiations, which gives
the desired result. It worth mentioning that our proof is in some sense
parallel to this calculation. For example, the number $N$ that appears
in the treatment of the second sum of Lemma \ref{rayslemma} is related
to the self-intersection numbers of the boundary divisors on the toric
surface given by the fans that correspond to the subgroups $S$.
\end{remark}

\begin{remark}
It may be interesting to analyze products of more than two $\tilde s$.
Every such product may be associated to a symbol
$$x_1^{r_1}\cdots x_n^{r_n}(a_1,\cdots,a_n)$$
where $a_i\in \ZZ/l\ZZ$. Then one expects to be able to develop a
generalization of the theory of Manin symbols, by introducing relations 
on these symbols that come from linear relations on the products.
The action of Hecke operators will then come from toric geometry, 
and will be related to subgroups of index $p$ in $\ZZ^{n}$ as in
\cite{toric}.
\end{remark}


\bibliographystyle{amsplain}
\bibliography{higher}

\end{document}